\documentclass[reqno]{amsart}

\usepackage{hyperref}
\usepackage{bbm}
\usepackage{amsfonts,amsmath,amsxtra,amsthm,amssymb,latexsym,enumerate}
\usepackage{ifpdf}
\usepackage{tikz}
\usetikzlibrary{patterns}
\usetikzlibrary{decorations.pathreplacing}
\usepackage{color}

\newcommand*{\mailto}[1]{\href{mailto:#1}{\nolinkurl{#1}}}
\newcommand{\arxiv}[1]{\href{http://arxiv.org/abs/#1}{arXiv:#1}}

\newcommand{\msc}[1]{\href{http://www.ams.org/msc/msc2010.html?t=&s=#1}{#1}}

\newtheorem{theorem}{Theorem}[section]
\newtheorem{corollary}[theorem]{Corollary}

\newtheorem{lemma}[theorem]{Lemma}
\newtheorem{hypothesis}{Hypothesis}[section]
\theoremstyle{definition}

\newtheorem{example}[theorem]{Example}
\newtheorem{remark}[theorem]{Remark}

\newcommand{\be}{\begin{equation}}
\newcommand{\ee}{\end{equation}}
\newcommand{\wh}{\widehat}
\newcommand{\id}{{\mathbbm 1}}

\numberwithin{equation}{section}

\newcommand{\floor}[1]{\lfloor#1 \rfloor}

\DeclareMathOperator{\dom}{dom}

\DeclareMathOperator{\loc}{loc}
\DeclareMathOperator{\grad}{grad}
\DeclareMathOperator{\Div}{div}

\newcommand\R{{\mathbb{R}}}
\newcommand\C{{\mathbb{C}}}
\newcommand\Z{{\mathbb{Z}}}

\newcommand{\gQ}{\mathfrak{Q}}

\newcommand\Ei{e_{\imath}}
\newcommand\Et{e_{\tau}}

\newcommand\Deg{{\rm{Deg}}}

\newcommand\cD{{\mathcal{D}}}

\newcommand\cI{{\mathcal{I}}}

\newcommand\cG{{\mathcal{G}}}
\newcommand\cE{{\mathcal{E}}}
\newcommand\cP{{\mathcal{P}}}
\newcommand\cV{{\mathcal{V}}}

\newcommand\bH{{\mathbf{H}}}
\newcommand\rh{{\mathbf{h}}}
\newcommand\rH{{\rm{H}}}
\newcommand\rJ{{\rm{J}}}
\newcommand\rI{{\rm{I}}}
\newcommand\E{{\rm{e}}}

\newcommand\Nr{{\rm{n}}}
\newcommand\I{{\rm{i}}}
\newcommand\rD{{\rm{d}}}

\newcommand\f{{\bf{f}}}

\def\wt#1{{{\widetilde #1} }}

\begin{document}

\title[Glazman--Povzner--Wienholtz on Graphs]{A Glazman--Povzner--Wienholtz Theorem on Graphs}

\author[A. Kostenko]{Aleksey Kostenko}
\address{Faculty of Mathematics and Physics\\ University of Ljubljana\\ Jadranska ul.\ 21\\ 1000 Ljubljana\\ Slovenia\\ and 
Faculty of Mathematics\\ University of Vienna\\
Oskar-Morgenstern-Platz 1\\ 1090 Vienna\\ Austria}
\email{\mailto{Aleksey.Kostenko@fmf.uni-lj.si}}

\author[M. Malamud]{Mark Malamud}
\address{RUDN University\\ Miklukho-Maklaya Str. 6\\ 117198 Moscow\\Russia}
\email{\mailto{malamud3m@gmail.com}}

\author[N. Nicolussi]{Noema Nicolussi}
\address{Centre de math\'ematiques Laurent Schwartz\\
\'Ecole Polytechnique\\ 91128 Palaiseau Cedex\\ France}  
\email{\mailto{noema.nicolussi@univie.ac.at}}

\thanks{{\it Research supported by the Austrian Science Fund (FWF) 
under Grants No.~I~4600~(A.K.) and J4497~(N.N.), and by the Slovenian Research Agency (ARRS) under Grant No.~N1-0137~(A.K.)}}
\thanks{\arxiv{2105.09931}}


\begin{abstract}
The Glazman--Povzner--Wienholtz theorem states that the completeness of a manifold, when combined with the semiboundedness of the Schr\"o\-din\-ger operator $-\Delta + q$ and suitable local regularity assumptions on $q$, guarantees its essential self-adjointness. Our aim is to extend this result to Schr\"odinger operators on graphs. We first obtain the corresponding theorem for Schr\"odinger operators on metric graphs, allowing in particular distributional potentials $q\in H^{-1}_{\rm loc}$. Moreover, we exploit recently discovered connections between Schr\"odinger operators on metric graphs and weighted graphs in order to prove a discrete version of the Glazman--Povzner--Wienholtz theorem.
\end{abstract}

\maketitle


\section{Introduction}

The problem of self-adjointness for Schr\"odinger-type operators
\begin{align}\label{eq:SchrGen}
\rH_q = -\Delta + q
\end{align}
is a classical topic having its origins in quantum mechanics and it is impossible to give even a very brief account on the existing literature (see, e.g., \cite{bms}, \cite{kato}, \cite{RSII}). Historically, the first work on the subject was done by H.~Weyl. Among many other results it was proved in~\cite{wey} that the 1d Schr\"odinger operator 
defined on the maximal domain in $L^2(\R)$ is self-adjoint if the real-valued potential $q\in C(\R)$ is bounded from below, $\inf_{\R}q(x) > -\infty$. The multidimensional version of this theorem turned out to be more subtle and it was proved later independently by T.~Carleman~\cite{car} and K.~Friedrichs~\cite{fri}.
 Seems, P.~Hartman~\cite{har48} (see also~\cite{rel51}) was the first to realize that in the case $n=1$ the semiboundedness of $q$ can be replaced by a much weaker condition -- it suffices to assume that the pre-minimal operator generated by the 1d Schr\"odinger operator in $L^2(\R)$ is bounded from below, that is, $\inf_{\R}q(x) > -\infty$ can be replaced by
\begin{align}
\inf_{f\in C_c^\infty(\R)\colon \|f\|_{L^2}=1} \int_\R |f'(x)|^2 + q(x)|f(x)|^2\rD x > -\infty.
\end{align}
The extension of this result to dimension $n\ge 2$ has an interesting history. To the best of our knowledge (see also \cite[Appendix~D.1]{bms} for further information), for Schr\"odinger operators in $\R^n$ the result was conjectured by I.M.~Glazman and proved by A.Ya.~Povzner in~1952 (see~\cite[\S~I.5]{pov53} and \cite[Theorem~35 on p.~58]{gla}). 
However, \cite{pov53} was published in Russian and was not widely known in the West until its English translation in 1967 (perhaps, it is not widely known even now). Being unaware of~\cite{pov53}, F.~Rellich in his invited address at the ICM in Amsterdam (1954) posed a multi-dimensional result as an open problem, which was solved a few years later by his student E.~Wienholtz~\cite{wie58}. Let us also emphasize that the approach suggested by Povzner in~\cite{pov53}, which is based on the relationship between essential self-adjointness and uniqueness of the Cauchy problem for the wave equation, was rediscovered later by P.R.~Chernoff~\cite{che}. 

There are (at least) two further possible directions for the above results. First of all, it makes sense to consider \eqref{eq:SchrGen} on Riemannian manifolds. On non-compact manifolds, the situation becomes nontrivial already when $q\equiv0$ since the geometry of a manifold comes into play. First results in this direction are due to M.P.~Gaffney~\cite{gaf}, \cite{gaf55} who noticed the importance of completeness of the manifold in question and the essential self-adjointness in this case was established later by W.~Roelcke \cite{roe} (see also \cite{che}, \cite{str})
\footnote{Very often the fact that completeness implies essential self-adjointness is attributed to M.P.~Gaffney, however, Gaffney in~\cite{gaf}, \cite{gaf55} actually proved the weaker property of Markovian uniqueness (that is, uniqueness of a Markovian extension). Namely, the maximal operator associated with $\Delta = \Div\cdot\grad$ is defined in ~\cite{gaf}, \cite{gaf55} as the product of the maximal operators $\Div$ and $\grad$ and hence it was assumed that functions from the maximal domain also belong to the Sobolev space $H^1$.}. Another direction is concerned with the local regularity assumptions on the potential $q$. The first significant improvement of the Carleman--Friedrichs result is due to T.~Kato \cite{kat74} (see also \cite[Theorem~X.28]{RSII}), who showed that $\rH_q$ is essentially self-adjoint if a semibounded potential $q$ is locally square integrable in $\R^n$, that is, one can replace the continuity assumption by $q\in L^2_{\loc}(\R^n)$. 
However, it turns out that it is impossible to replace $\inf_{\R^n}q(x) >-\infty$ by the semiboundedness of $\rH_q':=\rH_q|_{C_c^\infty}$ if $q$ is merely $L^2_{\loc}(\R^n)$ (indeed, if $n= 5$ and $q(x) = \alpha |x|^{-2}$ with $\alpha \in [-9/4,-5/4)$, then $\rH_q'$ is bounded from below, however, it is not essentially self-adjoint, see \cite[Example~D.1]{bms}, \cite[\S~X.4, p.185]{RSII}). In the particular case of dimension $n=1$, the situation seems to be rather complete in this respect. It was proved recently that Hartman's result remains true if $q$ is an $H^{-1}_{\loc}$ distribution (see \cite[Theorem 1 and Rem.~III.2]{akm10} and also \cite{hrmy}). 

After this very brief historical account let us turn to the main objects of this paper. Our main focus is on Schr\"odinger operators on graphs. We are interested in both metric graphs and weighted (discrete) graphs. Both objects have a venerable history and enjoy deep connections to several diverse branches of mathematics, placing them at the intersection of many subjects in mathematics and engineering. The key features of Laplacians on metric graphs, which are also widely known as {\em quantum graphs}, include their use as simplified models of complicated quantum systems (we only refer to \cite{bk13}, \cite{ekkst08} for further details and references).
The subject of {\em discrete Laplacians on graphs} is even wider and has been intensively studied from several perspectives (a partial overview of the immense literature can be found in \cite{cdv}, \cite{klwBook}, \cite{woe}). The analogs of the Gaffney and the Carleman--Friedrichs results in both contexts are rather well known. Indeed, in the discrete setting these results have been established in~\cite{hkmw13}, \cite{kl12}. In the context of metric graphs the situation is less transparent. On the one hand, the corresponding results seem to be a folklore and one can deduce them from, e.g., \cite{stu} (see \cite[Chap.~7.1]{kn21} for the details). On the other hand, we are aware of only two sources where the corresponding proofs can be found (see \cite{hae} and also \cite{ekmn}). As for the Glazman--Povzner--Wienholtz theorem on graphs, a discrete version of this result was proved only under additional restrictive assumptions on the geometry of the underlying graph, e.g., bounded geometry (see~\cite[Theorem~1.3]{mil11}, \cite[Theorem~6.1]{toha10})\footnote{After submission of this manuscript we learned from M.~Schmidt that the discrete version of the Glazman--Povzner--Wienholtz theorem without the bounded geometry assumption was proved in \cite[Theorem~2.16]{gks15}.}. In the context of metric graphs, such a result is not known to us even in the case of a continuous potential.  
In the present paper, we first establish the Glazman--Povzner--Wienholtz theorem on metric graphs allowing distributional potentials belonging to $H^{-1}_{\loc}$ (this, in particular, includes $\delta$-interactions on graphs). Moreover, using recently established connections between Laplacians on weighted graphs and on metric graphs, we employ this result to prove the discrete version of the Glazman--Povzner--Wienholtz theorem on arbitrary locally finite graphs, in particular, we remove the bounded geometry assumptions (Theorem~\ref{th:GPWdiscr}). Notice that the latter can be seen as a vast generalization of the classical self-adjointness test for Jacobi (tri-diagonal) matrices established by A.~Wouk~\cite{wouk} (see Remark~\ref{rem:wouk}).

Let us now outline the content. The first part of this paper gathers preliminary notions and facts. Namely, in Section~\ref{sec:prelim} we recall basic definitions and notation on graphs and metric graphs. Section~\ref{sec:GraphLapl} collects necessary information on weighted graph Laplacians and also on Schr\"odinger operators on graphs. First, we introduce the maximal and minimal weighted graph Laplacians on locally finite graphs. Then we recall the notion of an {\em intrinsic metric} on a weighted graph, which is necessary to state the Gaffney-type theorem (see Theorem~\ref{thm:GaffGraph}). Finally, we recall the definition of Schr\"odinger operators on graphs. The next Section~\ref{sec:LaplMetric} contains the information about Laplacians on metric graphs. We begin by introducing necessary function spaces on metric graphs and then continue with the definition of Kirchhoff Laplacians. It is important to underline that we are working in the setting of weighted metric graphs and do not make any assumptions regarding edge lengths. Next, we recall the notion of an intrinsic metric on a weighted metric graph (notice that the corresponding quadratic forms are strongly local Dirichlet forms) and complete this section by introducing Schr\"odinger operators on metric graphs. On the one hand, we prefer to work assuming some regularity on weights in order to streamline the exposition (see Hypothesis~\ref{hyp:munu}). On the other hand, we show how one can deal with low regularity potentials by allowing $q$ to be a real $H^{-1}_{\loc}$ distribution. 

The main results and their proofs are contained in Section~\ref{sec:GPWmetr} and Section~\ref{sec:GPWgraph}. Theorem~\ref{th:GPWcont} is a metric graph analog of the Glazman--Povzner--Wienholtz theorem. 
 Using recently established connections between Schr\"odinger operators on graphs and metric graphs (in order to make the exposition self-contained, we collect the necessary facts in Appendix~\ref{app:WGvsWMG}), Theorem~\ref{th:GPWcont} enables us to prove the Glazman--Povzner--Wienholtz theorem on weighted locally finite graphs (see Theorem~\ref{th:GPWdiscr}). We also apply this result to Jacobi matrices on graphs (Theorem~\ref{th:WoukGraph}), which can be seen as an extension to graphs of the classical self-adjointness test due to Wouk~\cite[Theorem~3(d)]{wouk} (see also \cite{Akh}).

\subsection*{Notation}
$\Z$, $\R$, $\C$ have their usual meaning; $\Z_{\ge 0} := \Z\cap [0,\infty)$.\\
$z^\ast$ denotes the complex conjugate of $z\in\C$.\\  
For a given set $S$, $\#S$ denotes its cardinality if $S$ is finite; otherwise we set $\#S=\infty$.\\
If it is not explicitly stated otherwise, we shall denote by $(x_n)$ a sequence $(x_n)_{n=0}^\infty$.

\section{Preliminaries} \label{sec:prelim}

\subsection{Graphs}\label{ss:II.01}

Let us first recall basic notions (we mainly follow the terminology in \cite{die}). Let $\cG_d = (\cV,\cE)$ be an undirected {\em graph}, that is, $\cV$ is a finite or countably infinite set of vertices and $\cE$ is a finite or countably infinite set of edges. 
Two vertices $u$, $v\in \mathcal{V}$ are called {\em neighbors} and we shall write $u\sim v$ if there is an edge $e_{u,v}\in \mathcal{E}$ connecting $u$ and $v$. For every $v\in \mathcal{V}$, we define $\cE_v$ as the set of edges incident to $v$. We stress that we allow {\em multigraphs}, that is, we allow {\em multiple edges} (two vertices can be joined by several edges) and {\em loops} (edges from one vertex to itself). Graphs without loops and multiple edges are called {\em simple}. 
Sometimes it is convenient to assign an {\em orientation} on $\cG_d$: to each edge $e\in\cE$ one assigns the pair $(\Ei,\Et)$ of its {\em initial} $\Ei$ and {\em terminal} $\Et$ vertices. We shall denote the corresponding oriented graph by $\vec{\cG}_d = (\cV,\vec{\cE})$, where $\vec{\cE}$ denotes the set of oriented edges. Notice that for an oriented loop we do distinguish between its initial and terminal vertices. Next, for every vertex $v\in\cV$, set 
\begin{align}\label{eq:Ev_pm}
{\cE}^+_v & = \big\{(\Ei,\Et) \in \vec{\cE}\,|\, \Ei = v\big\}, & {\cE}_v^- & = \big\{(\Ei,\Et) \in \vec{\cE} \,|\,  \Et = v\big\},
\end{align}
and let $\vec{\cE}_v$ be the disjoint union of outgoing $\cE_v^+$ and incoming $\cE_v^-$ edges,
\begin{align}\label{eq:vecEv}
\vec{\cE}_v & := {\cE}_v^+ \sqcup {\cE}_v^- = \vec{\cE}_v^+ \cup \vec{\cE}_v^-, &  \vec{\cE}_v^\pm & := \big\{(\pm,e)\,|\, e\in \cE_v^\pm\big\}.
\end{align}
We shall denote the elements of $\vec{\cE}_v$ by $\vec{e}$ and also introduce the orientation function $\pi_v\colon \vec{\cE}_v\to \{-1,1\}$  by setting $\pi_v(\vec{e}) = \pm 1$ for $\vec{e}\in \vec{\cE}_v^\pm$. 
 
The {\em (combinatorial) degree} of $v\in\cV$ is 
\begin{align}\label{eq:combdeg}
\deg(v):=  \#(\vec{\cE}_v ) = \#(\vec{\cE}_v^+ ) + \#(\vec{\cE}_v^- )  = \#(\cE_v ) + \#\{e\in\cE_v|\, e\ \text{is a loop}\}.
\end{align}
Notice that if $\cE_v$ has no loops, then $\deg(v) = \#(\cE_v)$. The graph $\cG_d$ is called {\em locally finite}  if $\deg(v)<\infty$ for all $v\in\cV$. 

 A sequence of (unoriented) edges $\cP = (e_{v_0, v_1}, e_{v_1, v_2}, \dots, e_{v_{n-1}, v_n})$ is called a {\em path} of (combinatorial) length $n\in \Z_{\ge 0}\cup \{\infty\}$.  Notice that for simple graphs each path $\cP$ can be identified with its sequence of vertices $\cP = (v_k)_{k=0}^n$.
 A graph $\cG_d$ is called {\em connected} if for any  two vertices there is a path connecting them. 
 
We shall always make the following assumptions on the geometry of $\cG_d$:

\begin{hypothesis}\label{hyp:graph01}
$\cG_d$ is connected and locally finite. 
\end{hypothesis}

\subsection{Metric graphs}\label{ss:II.02}

Assigning each edge $e\in\cE$ a finite length $|e| \in (0,\infty)$, we can naturally associate with $(\cG_d,|\cdot|) = (\cV,\cE,|\cdot|)$ a metric space $\cG$: first, we identify each edge $e \in \cE$ with the copy of the interval $\cI_e = [0, |e|]$, which also assigns an orientation on $\cE$ upon identification of $\Ei$ and $\Et$ with the left, respectively, right endpoint of $\cI_e$. The topological space $\cG$ is then obtained by ``glueing together" the ends of edges corresponding to the same vertex $v$ (in the sense of a topological quotient, see, e.g., \cite[Chap.~3.2.2]{bbi}). 
The topology on $\cG$ is metrizable by the {\em length metric} $\varrho_0$ --- the distance between two points $x,y \in\cG$ is defined as the arc length of the ``shortest path" connecting them (such a path does not necessarily exist and one needs to take the infimum over all paths connecting $x$ and $y$). 

A \emph{metric graph} is a (locally compact) metric space $\cG$ arising from the above construction for some collection $(\cG_d, |\cdot|) =(\cV, \cE, |\cdot|)$. More specifically, $\cG$ is then called the \emph{metric realisation} of $(\cG_d, |\cdot|)$, and a pair $(\cG_d, |\cdot|)$ whose metric realization coincides with $\cG$ is called a \emph{model} of $\cG$. For a thorough discussion of metric graphs as topological and metric spaces we refer to \cite[Chap.~I]{hae}. 

\begin{remark}\label{rem:II.mr=lengthspace}
Let us stress that a metric graph $\cG$ equipped with the length metric $\varrho_0$ (or with any other path metric) is a {\em length space} (see \cite[Chap.~2.1]{bbi} for definitions and further details). In particular, taking into account our Hypothesis~\ref{hyp:graph01}, the Hopf--Rinow theorem holds true on $\cG$ (cf.~\cite[Theorem~2.5.28]{bbi}), which relates completeness with both geodesic completeness and bounded compactness. Let us also mention that complete, locally compact length spaces are {\em geodesic}, that is, every two points can be connected by a shortest path \cite[Theorem~2.5.23]{bbi}.
\end{remark}

Clearly, different models may give rise to the same metric graph. Moreover,  any metric graph has infinitely many models (e.g., they can be constructed by subdividing edges using vertices of degree two). 
A model $(\cV,\cE, |\cdot|)$ is called {\em simple} if the corresponding graph $(\cV,\cE)$ is simple. In particular, every locally finite metric graph has a simple model and hence this indicates that restricting to simple graphs, that is, assuming in addition to Hypothesis~\ref{hyp:graph01} that $\cG_d$ has no loops or multiple edges, would not be a restriction at all when dealing with metric graphs. 

\begin{remark}\label{rem:Models}
In most parts of our paper, we will consider a metric graph together with a fixed choice of its model. In this situation, we will usually be slightly imprecise and do not distinguish between these two objects. In particular, we will denote both objects by the same letter $\cG$ and write $\cG = (\cV, \cE, |\cdot|)$ or $\cG= (\cG_d,|\cdot|)$. 
\end{remark}

\begin{remark}[Metric graph as a 1d manifold with singularities]\label{rem:MGasM1}
Sometimes it is useful to consider metric graphs as one-dimensional manifolds with singularities.
Since every point $x\in\cG$ has a neighborhood isomorphic to a star shaped set
\begin{align} \label{eq:star}
\cE(\deg(x),r_x) := \big\{z= r\E^{2\pi \I k/\deg(x)}|\, r\in [0,r_x),\ k=1,\dots,\deg(x) \big\}\subset \C,
\end{align}
 one may introduce the set of {\em tangential directions} $T_x(\cG) $ at $x$ as the set of unit vectors $\E^{2\pi \I k/\deg(x)}$, $k=1,\dots, \deg(x)$. 
Then all vertices $v\in\cV$ with $\deg(v)\ge 3$ are considered as {\em branching points/singularities}
and vertices $v\in\cV$ with $\deg(v)=1$ as {\em boundary points}. 
Notice that for every vertex $v\in\cV$ the set of tangential directions $T_v(\cG)$ can be identified with $\vec{\cE}_v$. 
 If there are no loop edges at the vertex $v \in \cV$, then $T_v(\cG)$ is 
identified with $\cE_v$ in this way. 
\end{remark}

\section{Discrete Laplacians and Schr\"odinger-type operators on graphs}\label{sec:GraphLapl}  

\subsection{Graph Laplacians} \label{ss:III.01}

Let $\cV$ be a countable set. A function $m\colon \cV\to (0,\infty)$ defines a measure of full support on $\cV$ in an obvious way. A pair $(\cV,m)$ is called a {\em discrete measure space}. The set of square summable functions 
\[
\ell^2(\cV;m) = \Big\{ f\in C(\cV)\,|\,\, \|f\|^2_{\ell^2(\cV;m)}:= \sum_{v\in\cV} |f(v)|^2 m(v) <\infty \Big\}
\]
has a natural Hilbert space structure. Here $C(\cV)$ is the set of complex-valued functions on a countable set $\cV$. 

Suppose $b\colon \cV\times\cV \to [0,\infty)$ satisfies the following conditions:
\begin{itemize}
\item[(i)] {\em symmetry}: $b(u,v) = b(v,u)$ for each pair $(u,v)\in \cV\times\cV$,
\item[(ii)] {\em vanishing diagonal}: $b(v,v) = 0$ for all $v\in\cV$,
\item[(iii)] {\em locally finite}:   $\#\{ u\in\cV\,|\, b(u,v)\neq 0\} < \infty$  for all $v\in\cV$\footnote{In fact, using the form approach one can considerably relax this condition by replacing it  with the {\em local summability}: $\sum_{v\in\cV} b(u,v) <\infty$ for all $u\in\cV$.}.
\item[(iv)] {\em connected}: for any $u,v\in \cV$ there is a finite collection $(v_k)_{k=0}^n \subset \cV$ such that $u=v_0$, $v=v_n$ and $b(v_{k-1},v_k)>0$ for all $k\in \{1,\dots,n\}$.
\end{itemize}
Following \cite{kl12}, $b$ is called a {\em (weighted) graph} over $\cV$ or over $(\cV,m)$ if in addition  a measure $m$ of full support on $\cV$ is given ($b$ is also called an {\em edge weight}).  To simplify notation, we shall denote a graph $b$  over $(\cV,m)$ by $(\cV,m;b)$. 

\begin{remark}\label{rem:simplevsmult}
Let us quickly explain how the above notion is related to the previous section. To any graph $b$ over $\cV$, we can naturally associate a simple combinatorial graph $\cG_b$. Namely, the vertex set of $\cG_b$ is $\cV$ and its edge set $\cE_b$ is defined by calling two vertices $u,v\in\cV$ neighbors, $u\sim v$, exactly when $b(u,v)>0$. Clearly, $\cG_b = (\cV,\cE_b)$ is an undirected graph in the sense of Section~\ref{ss:II.01}. Let us stress, however, that the constructed {\em graph $\cG_b$ is always simple}. 
\end{remark}

The {\em (formal) Laplacian} $L = L_{m,b}$ associated to a graph $b$ over $(\cV,m)$ is given by
\begin{align}\label{eq:LaplDiscr}
(L f)(v) = \frac{1}{m(v)}\sum_{u\in\cV} b(v,u) (f(v) - f(u)).
\end{align} 
It acts on functions $f\in C(\cV)$ and this naturally leads to the {\em maximal} Laplacian $\rh$ in $\ell^2(\cV;m)$ defined by 
\begin{align}\label{eq:LaplDiscrMax}
\rh & = L\upharpoonright\dom(\rh), & \dom(\rh) & = \{f\in \ell^2(\cV;m)\,|\,  Lf\in\ell^2(\cV;m)\}.
\end{align}
This operator is closed, however, if $\cV$ is infinite, it is not symmetric in general (cf. \cite[Theorem~6]{kl12}). 
Taking into account that $b$ is locally finite, it is straightforward to verify that $C_c(\cV)\subseteq \dom(\rh)$. Therefore, we can introduce the {\em minimal} Laplacian $\rh^0$ as the closure in $\ell^2(\cV;m)$ of the {\em pre-minimal} Laplacian
\begin{align}\label{eq:pmLaplDiscr}
\rh' = L\upharpoonright\dom(\rh'),\qquad \dom(\rh')= C_c(\cV).
\end{align}
Then $\rh'\subseteq \rh^0\subseteq \rh$ and $(\rh')^\ast = (\rh^0)^\ast  = \rh$.  
The following fact is rather well known (see, e.g., \cite[Lemma~1]{dav}, \cite[Theorem~11]{kl10}, \cite[Rem.~1]{susy}). 

\begin{lemma}\label{lem:boundeddiscr} 
The Laplacian $L = L_{m,b}$ is bounded on $\ell^2(\cV,m)$ if and only if the weighted degree function
  $\Deg\colon \cV \to [0,\infty)$ given by 
\begin{align}\label{eq:WDegDef}
\Deg\colon v\mapsto \frac{1}{m(v)} \sum_{u\in \cV} b(u,v)
\end{align}
is bounded on $\cV$. In this case $\rh^0 = \rh$ and $\|\Deg\|_{\infty}\le \|\rh\|_{\ell^2(\cV;m)} \le 2\|\Deg\|_{\infty}$. 
\end{lemma}

\subsection{Intrinsic metrics and Gaffney-type theorems on graphs}\label{ss:III.02}

Let $b$ be a connected graph over $(\cV, m)$. 
A symmetric function $p\colon \cV \times \cV \to [0,\infty)$ such that $p(u,v)>0$ exactly when $b(u,v) >0$ is called a \emph{weight function} for $(\cV,m;b)$. Every weight function $p$ generates a {\em path metric} $\varrho_p$ on $\cV$ with respect to the graph $b$ via
\begin{align}\label{eq:pathmetric}
\varrho_p(u,v) := 
\inf_{\cP= (v_0,\dots,v_n)\colon u=v_0,\ v=v_n}\sum_{k} p(v_{k-1},v_k). 
\end{align}
Here the infimum is taken over all paths in $(\cV,m;b)$ connecting $u$ and $v$.
Since $b$ is locally finite, $\varrho_p(u,v) > 0$ whenever $u \neq v$.

Following \cite{flw14} (see also \cite{kel15}), a metric $\varrho$ is called {\em intrinsic} w.r.t. $(\cV,m;b)$ if 
\begin{align}\label{eq:intrinsicdef}
 \sum_{u\in\cV}b (u,v)\varrho(u,v)^2\le m (v)
\end{align}
holds for all  $v\in\cV$. Similarly, 
a weight function $p\colon \cV \times \cV \to [0,\infty)$ is called an \emph{intrinsic weight} for $(\cV,m;b)$ if 
\begin{align}\label{eq:intrweigdef}
	 \sum_{u\in\cV}b (u,v) p(u,v)^2\le m (v),\qquad v\in\cV.
\end{align}
If $p$ is an intrinsic weight, then the associated path metric $\varrho_p$ is called {\em strongly intrinsic} (it is obviously intrinsic in the sense of \eqref{eq:intrinsicdef}).

\begin{remark} \label{rem:intrinsic_existence}
For any given locally finite graph $(\cV,m;b)$ an intrinsic metric always exists (see \cite[Example~2.1]{hkmw13}, \cite{kel15}, \cite{dVTHT} and also \cite[Chap.~6]{kn21}).
\end{remark}

It is well known in the context of manifolds that completeness of a Riemannian manifold implies self-adjointness of the corresponding Laplace--Beltrami operator. The analog of this result in the graph setting also holds true:

\begin{theorem}[Gaffney's theorem on graphs]\label{thm:GaffGraph}
Let $b$ be a locally finite, connected graph over $(\cV,m)$ and let $\rh$ be the maximal Laplace operator in $\ell^2(\cV;m)$. Then $\rh$ is self-adjoint if at least one of the following two conditions is satisfied:
\begin{itemize}
\item[(i)]
$(\cV,\varrho)$ is complete as a metric space for some intrinsic path metric  $\varrho$,
\item[(ii)]
$(\cV,\varrho_m)$ is complete, where $\varrho_m$ is the metric given by
\begin{align}\label{eq:rho_m}
\varrho_m(u,v) = \inf_{\cP = (v_k)_k}\sum_{k} m(v_{k}),
\end{align}
where the infimum is taken over all paths in $(\cV,m;b)$ connecting $u$ and $v$.
\end{itemize}
\end{theorem}

\begin{remark}
The sufficiency of (i) for the self-adjointness was established in~\cite{hkmw13}; (ii) was proved in~\cite{kl12}. Let us stress that~\cite{hkmw13} and~\cite{kl12} also obtain the corresponding analogs for non-locally finite edge weights $b$, assuming only the {\em local summability condition}. 
Other necessary and sufficient self-adjointness conditions can be found in~\cite{dVTHT}, \cite{dVTHT3}, \cite{mil11}, \cite{toha10}.
\end{remark}

\begin{example}[Combinatorial Laplacian]\label{ex:CombLapl}
Let $m\equiv \id$ on $\cV$ and $b$ be the adjacency matrix of a (simple) graph, $b\colon \cV\times\cV\to \{0,1\}$. The corresponding Laplacian~\eqref{eq:LaplDiscr}
\begin{align}\label{eq:LaplComb}
(L_{\rm comb} f)(v) = \sum_{u\sim v} f(v) - f(u) = \deg(v)f(v) - \sum_{u\sim v}f(u)
\end{align} 
is called the {\em combinatorial Laplacian}. The maximal operator $\rh_{\rm comb}$ is bounded on $\ell^2(\cV)$ exactly when $\cG_b$ has bounded geometry, $\sup_{v\in\cV} \deg(v)<\infty$. However, its self-adjointness, proved independently in \cite{jor}, \cite{weber10}, \cite{woj}, immediately follows from Theorem~\ref{thm:GaffGraph}(ii). Notice that Theorem~\ref{thm:GaffGraph}(i) enables us to conclude the self-adjointness of $\rh_{\rm comb}$ only if the degree function does not grow too fast at ``infinity". 
\end{example}

\subsection{Schr\"odinger-type operators on graphs}\label{ss:III.03}

Let $\alpha\colon\cV\to \R$ be a real-valued function on $\cV$. Following the considerations of Section \ref{ss:III.01}, we can associate in $\ell^2(\cV;m)$ the maximal and minimal operators with the Schr\"odinger-type expression
\begin{align}\label{eq:SchrDiscr}
(L_\alpha f)(v) = (L_{m,b,\alpha} f)(v):= \frac{1}{m(v)}\Big(\sum_{u\in\cV} b(v,u) (f(v) - f(u)) + \alpha(v) f(v)\Big).
\end{align} 
Namely, $L_\alpha$ is well-defined on $C(\cV)$ and then the maximal Schr\"odinger operator $\rh_\alpha$ is defined as a restriction of $L_\alpha$ to $\ell^2(\cV;m)$:
\begin{align}\label{eq:SchrDiscrMax}
\rh_\alpha & = L_\alpha\upharpoonright\dom(\rh_\alpha), & \dom(\rh_\alpha) & = \{f\in \ell^2(\cV;m)\,|\, L_\alpha f\in\ell^2(\cV;m)\}.
\end{align}
The pre-minimal $\rh_\alpha'$ and minimal $\rh_\alpha^0$ operators are defined analogously. It turns out that Theorem \ref{thm:GaffGraph} extends to Schr\"odinger operators on graphs, however, under the additional assumption that the potential $\alpha$ is bounded from below (see \cite{hkmw13}, \cite{kl12}), that is,
\begin{align}
\inf_{\cV} \frac{\alpha(v)}{m(v)} > -\infty.
\end{align}
The latter can be seen as the graph analog of the Carleman--Friedrichs theorem for Schr\"odinger operators in $\R^n$.

\section{Laplacians and Schr\"odinger operators on metric graphs}\label{sec:LaplMetric}

\subsection{Function spaces on metric graphs}\label{ss:IV.01}  
Let $\cG$ be a metric graph together with a fixed model $(\cV,\cE,|\cdot|)$. Identifying every edge $e\in\cE$ with a copy of $\cI_e = [0,|e|]$, we can introduce Lebesgue  spaces on edges and also on $\cG$. 
First of all, let $\mu\colon \cG\to (0,\infty)$ be a measurable function  on $\cG$. 
Considering the weight $\mu$ as a density, we can associate with it a measure $\mu$ on $\cG$, which is given by $\mu(\rD x) = \mu(x_e)\rD x_e$ on every edge $e\in\cE$. Thus, we can define the Hilbert space $L^2(\cG;\mu)$ of functions $f\colon \cG\to \C$ which are square integrable w.r.t the measure $\mu$ on $\cG$. Similarly, one defines the Banach spaces $L^p(\cG;\mu)$ for any $p\in [1,\infty]$. In fact, if $p\in [1,\infty)$, then 
$L^p(\cG;\mu)$ can be seen as the edgewise direct sum of $L^p$ spaces 
\[
L^p(\cG;\mu) \cong  \Big\{f=(f_e)_{e\in\cE}\big|\, f_e\in L^p(e;\mu),\ \sum_{e\in\cE}\|f_e\|^p_{L^p(e;\mu)}<\infty\Big\},
\]
where 
\[
\|f_e\|^p_{L^p(e;\mu)} = \int_{e}|f_e(x_e)|^p\mu(\rD x_e) =  \int_{e}|f_e(x_e)|^p\,\mu(x_e)\rD x_e,
\]
that is, $L^p(e;\mu)$ stands for the usual $L^p$ space upon identifying $e$ with $\cI_e$ and $\mu$ with the measure $\mu(x_e)\rD x_e$ on $\cI_e$. If $\mu = 1$ on $e$, then we shall simply write $L^p(e)$. 
Next, the subspace of compactly supported $L^p$ functions will be denoted by
\begin{equation*}
	L^p_c(\cG;\mu) = \big\{f \in L^p(\cG;\mu)\,| \; f|_e \not\equiv 0 \text{ only on finitely many edges } e \in \cE\big\}.
\end{equation*}
The space  $L^p_{\loc}(\cG;\mu)$  of locally $L^p$ functions consists of all measurable functions $f$ such that $fg\in L^p_c(\cG;\mu)$ for all $g\in C_c(\cG)$. Notice that the weight $\mu$ is edgewise $L^1$ if and only if $\mu \in L^1_{\loc}(\cG)$.
 
\subsection{Kirchhoff Laplacians}\label{ss:IV.02}

Again, let $\cG$ be a metric graph. 
 Suppose we are also given two positive measurable functions 
\begin{align}
\mu \colon & \cG \to  (0,\infty), & 
\nu \colon & \cG \to  (0,\infty).
\end{align}
In the following, if not stated otherwise, we shall impose the following regularity assumption on the weights:

\begin{hypothesis}\label{hyp:munu}
The weights $\mu$ and $\nu$ as well as their reciprocals $1/\mu$ and $1/\nu$ are edgewise continuously differentiable\footnote{We shall say that a function $f\colon \cG\to \C$ is {\em edgewise continuously differentiable} on $\cG$ if there is a model $(\cV,\cE,|\cdot|)$ of $\cG$ such that the restriction of $f$ to the interior of each edge $e$ admits a continuation to a continuously differentiable function on $e$.} on $\cG$.
\end{hypothesis}

\begin{remark}\label{rem:munu1}
This regularity assumption on edge weights is not optimal and one can,  of course, consider locally integrable weights, i.e., assume that $\mu,1/\nu\in L^1_{\loc}(\cG)$. However, we decided to restrict to the class of edgewise $C^1$ weights in order to streamline the exposition. In the following, at the corresponding places  we shall indicate, omitting the details, how to modify the arguments in order to achieve weaker regularity assumptions.
\end{remark}

Fix a model of $\cG$ and for every $e\in\cE$ consider the maximal operator $\rH_{e,\max}$ defined in $L^2(e;\mu)$ by
\begin{align}\label{eq:Hemax}
 \rH_{e,\max}f & = \tau_e f,\qquad \tau_e = -\frac{1}{\mu(x_e)}\frac{\rD}{\rD x_e}\nu(x_e)\frac{\rD}{\rD x_e}, \\
 \dom(\rH_{e,\max}) & = \big\{f\in L^2(e;\mu)\,|\, f,\ \nu f'\in AC(e),\ \tau_e f\in  L^2(e;\mu) \big\}.
\end{align}
Taking into account our Hypothesis~\ref{hyp:munu}, it is straightforward to show that in fact $\dom(\rH_{e,\max})$ coincides with the Sobolev space $H^2(e)$ algebraically and topologically. The maximal operator on $\cG$ is then defined in $L^2(\cG;\mu)$ as
\begin{align}\label{eq:Hmax}
\bH_{\max}  = \bigoplus_{e\in \cE} \rH_{e,\max}. 
\end{align}
Clearly, for each $f \in \dom(\bH_{\max})$ the following quantities
 \begin{align}\label{eq:tr_fe}
 f(\Ei) & := \lim_{x_e \to \Ei} f (x_e), & f (\Et) & := \lim_{x_e \to  \Et} f (x_e),
 \end{align}
 and the normal derivatives 
 \begin{align}\label{eq:tr_fe'}
 \partial_{\vec{e}} f (\Ei) & := \lim_{x_e \to \Ei} \frac{f (x_e) - f(\Ei)}{|x_e - \Ei|}, & 
 \partial_{\vec{e}} f (\Et) & := \lim_{x_e \to \Et} \frac{f (x_e) - f(\Et)}{|x_e - \Et|},
 \end{align}
are well defined for all $e\in\cE$. We also need the following notation
\begin{align}
f_{\vec{e}}(v) & := \begin{cases} f(\Ei), & \vec{e}\in \vec{\cE}_v^+, \\ f(\Et), & \vec{e} \in \vec{\cE}_v^-, \end{cases}
& \partial_{\vec{e}} f(v) & := \begin{cases} \partial f (\Ei), & \vec{e}\in \vec{\cE}_v^+, \\ \partial f (\Et), & \vec{e}\in \vec{\cE}_v^-, \end{cases}
\end{align}
for every $v\in\cV$ and $\vec{e}\in\vec{\cE}_v$. In the case of a loopless graph, the above notation simplifies since we can identify $\vec{\cE}_v$ with $\cE_v$ for all $v\in\cV$.  

Now, in order to reflect the underlying graph structure, we impose at each 
vertex $v\in\cV$ the following boundary conditions
\begin{align}\label{eq:vert-alpha}
\begin{cases} f\ \text{is continuous at}\ v,\\[1mm] 
\sum\limits_{\vec{e}\in \vec{\cE}_v} \nu_{\vec{e}}(v)\partial_{\vec{e}} f(v) = \alpha(v)f(v), \end{cases} 
\end{align} 
where $\alpha(v)\in \R\cup\{\infty\}$, and $\alpha(v) = \infty$ should be understood as the Dirichlet boundary condition at $v$. These vertex conditions are interpreted as {\em $\delta$-interactions of strength} $\alpha$ (see Remark~\ref{rem:LaplDeltaII}). 

\begin{remark}
Replacing Hypothesis~\ref{hyp:munu} by the $L^1_{\loc}$ condition on $\mu$ and $1/\nu$, one then simply needs to modify \eqref{eq:tr_fe'} by considering edgewise quasi-derivatives $\nu f'$ since on each edge $e$ the function $\nu f'$ is absolutely continuous if $f\in \dom(\rH_{e,\max})$.
\end{remark}

To motivate our definition, for edgewise locally absolutely continuous functions on $\cG$, let us denote by $\nabla$ the edgewise first derivative,
\begin{align}\label{eq:nabla}
\nabla\colon f\mapsto f'
\end{align}
Next, consider $\nabla$ as the differentiation operator on $\cG$ acting on functions which are edgewise locally absolutely continuous and also continuous at the vertices.
Notice that when considering $\nabla$ as an operator acting from $L^2(\cG;\mu)$ to $L^2(\cG;\nu)$, its formal adjoint $\nabla^\dagger$ acting from $L^2(\cG;\nu)$ to $L^2(\cG;\mu)$ acts edgewise as
\begin{align}
\nabla^\dagger \colon f\mapsto -\frac{1}{\mu} (\nu f)'.
\end{align}
Thus, the weighted Laplacian $\Delta$ acting in $L^2(\cG;\mu)$, written in the divergence form 
\begin{align}\label{eq:distrLapl}
\Delta \colon f\mapsto  -\nabla^\dagger (\nabla f),
\end{align}
acts edgewise as the following divergence form Sturm--Liouville operator
\begin{align}\label{eq:LaplMetrG}
\Delta\colon f\mapsto \frac{1}{\mu}(\nu f')'.
\end{align}
The continuity assumption imposed on $f$ results for $\Delta$ in a one-parameter family of symmetric boundary conditions \eqref{eq:vert-alpha}. With the Laplacian $\Delta$ acting on $\cG$ we shall always associate  
the {\em Kirchhoff} vertex conditions:
\begin{align}\label{eq:kirchhoff}
\begin{cases} f\ \text{is continuous at}\ v,\\[1mm] 
\sum\limits_{\vec{e}\in \vec{\cE}_v} \nu_{\vec{e}}(v)\partial_{\vec{e}} f(v) = 0, \end{cases}\qquad v\in\cV,
\end{align} 
that is, conditions~\eqref{eq:vert-alpha} with $\alpha(v) = 0$ for all $v\in\cV$. 
In particular, imposing these boundary conditions on the maximal domain yields the \emph{(maximal) Kirchhoff Laplacian}:
\begin{align}\label{eq:H}
\begin{split}
	\bH  =  -\Delta\upharpoonright &{\dom(\bH)},\\
	& \dom(\bH ) = \{f\in \dom(\bH_{\max})\,|\, f\ \text{satisfies}~\eqref{eq:kirchhoff}\ \text{on}\ \cV\}.
\end{split}
\end{align}

The maximal operator with boundary conditions \eqref{eq:vert-alpha} will be denoted by $\bH_\alpha$:
\begin{align}\label{eq:Halpha}
\begin{split}
	\bH_{\alpha}   =  -\Delta\upharpoonright &{\dom(\bH_{\alpha})},\\ 
	 & \dom(\bH_{\alpha} )  = \{f\in \dom(\bH_{\max})\,|\, f\ \text{satisfies}~\eqref{eq:vert-alpha}\ \text{on}\ \cV\}.
\end{split}
\end{align} 
 
\subsection{Energy forms} \label{ss:IV.03}

Restricting further to compactly supported functions we end up with the pre-minimal operator
\begin{align}\label{eq:Halpha0}
	\bH_{\alpha}'  =  -\Delta\upharpoonright {\dom(\bH_{\alpha}')},\qquad 
	 \dom(\bH_{\alpha}')  = \dom(\bH_{\alpha}) \cap C_c(\cG).
\end{align}
Integrating by parts one obtains
\begin{align}\label{eq:QFalpha}
	\langle\bH_\alpha' f, f\rangle_{L^2} = \int_\cG |\nabla f(x)|^2 \, \nu(\rD x) + \sum_{v\in\cV}\alpha(v)|f(v)|^2 =: \gQ_\alpha[f] 
\end{align}
for all $f \in \dom(\bH_\alpha')$, 
which implies that $\bH_\alpha'$ is a symmetric operator in $L^2(\cG;\mu)$. We define $\bH_\alpha^0$ as the closure of $\bH_{\alpha}'$ in $L^2(\cG;\mu)$. It is standard to show that
\begin{align}\label{eq:H0*=HA}
		(\bH_\alpha')^\ast = \bH_\alpha.
\end{align}
In particular, the equality $\bH_\alpha^0 = \bH_\alpha$ holds if and only if $\bH_\alpha$ is self-adjoint (or, equivalently, $\bH_\alpha'$ is essentially self-adjoint). 

If $\alpha\equiv 0$, that is, the case of the Kirchhoff Laplacian, the corresponding minimal and pre-minimal operators will be denoted by $\bH^0$ and $\bH'$. The energy form in this case is simply the $L^2$ norm of the gradient
\begin{align}\label{eq:integrationbp}
	\gQ[f] = \int_\cG |\nabla f(x)|^2\nu(\rD x),
\end{align}
which implies that both $\bH'$ and $\bH^0$ are nonnegative symmetric operators. With this form we shall associate two spaces: the first Sobolev space $H^1(\cG) = H^1(\cG;\mu,\nu)$ is defined as the subspace of $L^2(\cG;\mu)$ consisting of continuous functions, which are edgewise absolutely continuous and have finite energy $\gQ[f] < \infty$. Equipping $H^1$ with the standard graph norm turns it into a Hilbert space. Also, we define the space $H^1_0(\cG) = H^1_0(\cG;\mu,\nu)$ as the closure of compactly supported $H^1$ functions in the $H^1$ norm,
\[
H^1_0 = H^1_0(\cG;\mu,\nu):= \overline{H^1_c(\cG)}^{\|\cdot\|_{H^1(\cG;\mu,\nu)}},
\]
where $H^1_c(\cG) := H^1(\cG)\cap C_c(\cG)$. 
Restricting $\gQ$ to these spaces, we end up with two closed forms in $L^2(\cG;\mu)$:
\begin{align}
\gQ_D & = \gQ\upharpoonright{H^1_0}, & \gQ_N & = \gQ\upharpoonright{H^1}.
\end{align}
According to the representation theorem, they give rise to two self-adjoint nonnegative operators $\bH_D$ and $\bH_N$ in $L^2(\cG;\mu)$, the Dirichlet and Neumann Laplacians, respectively. Notice also that $\bH_D$ coincides with the Friedrichs extension of $\bH'$:
\[
\dom(\bH_D) = \dom(\bH)\cap H^1_0(\cG).
\]

\subsection{Intrinsic metrics on metric graphs}\label{ss:VI.04}

We define the intrinsic metric $\varrho$ of a weighted metric graph $(\cG, \mu, \nu)$ as the (maximal) intrinsic metric of its Dirichlet Laplacian $\bH_D$ (in particular, note that $\gQ_D$ is a strongly local, regular Dirichlet form). By \cite[eq.~(1.3)]{stu} (see also \cite[Theorem~6.1]{flw14}), $\varrho_{\rm intr}$ is given by
\[
	\varrho_{\rm intr}(x,y) = \sup \big \{ f(x) - f(y)\, | \, f \in \wh{\cD}_{\loc} \big \}, \qquad x, y \in \cG, 
\]
where the function space $\wh{\cD}_{\loc}$ is defined as
\[
\wh{\cD}_{\loc} = \big\{ f \in H^1_{\loc}(\cG)\, \big | \ \nu(x) |\nabla f (x) |^2 \le \mu(x)\ \ \text{for a.e.}\ x\in\cG \big\}.
\]
Here $H^1_{\loc}(\cG) = \{f\in C(\cG)|\, fg\in H^1(\cG)\ \text{for all}\ g\in H^1_c(\cG)\}$. 
It turns out that the intrinsic metric $\varrho$ admits a rather explicit description. First of all, the above suggests to define the intrinsic weight $\eta\colon \cG \to (0,\infty)$,
\begin{align}
\eta = \eta_{\mu,\nu} := \sqrt{ \frac{\mu }{ \nu}}\quad \text{on}\ \ \cG.
\end{align}
This weight gives rise to a new measure on $\cG$ whose density w.r.t. the Lebesgue measure is exactly $\eta$ (as in the case of the edge weights on a metric graph, we abuse the notation and denote with $\eta$ both the edge weight and the corresponding measure). Due to our assumptions on both $\mu$ and $\nu$, the weight $\eta$ is positive a.e. on $\cG$ and, moreover, edgewise $C^1$ on $\cG$. 
Next, take any path $\cP$ in $\cG$ (taking into account the local structure of $\cG$, it suffices to assume that $\cP$ is an image of a continuous, piecewise injective map $\gamma\colon I\to \cG$, where $I\subseteq \R$ is an interval). The \emph{(weighted) length} of a path $\cP \subset \cG$ is defined as
\begin{align}
	|\cP|_{\eta}  := \int_{\cP} \eta(\rD x) = \int_{\cP} \sqrt{\frac{\mu}{\nu}}\,\rD x.  
\end{align}

\begin{lemma}\label{lem:IntrMetWMG}
Let $\cG$ be a metric graph and let $\mu$, $\nu$ be two positive weights such that $\mu$, $1/\nu \in L^1_{\loc}(\cG)$. 
Then the metric $\varrho_\eta$ defined by 
\begin{equation} \label{eq:compute_metric}
\varrho_\eta(x,y) :=\inf_\cP |\cP|_\eta = \inf_\cP \int_{\cP} \eta(\rD x), \qquad x, y \in \cG,
\end{equation} 
where the infimum is taken over all paths $\cP$ from $x$ to $y$, coincides with the intrinsic metric on $(\cG, \mu, \nu)$ (w.r.t. $\gQ_D$), that is, $\varrho_{\rm intr} = \varrho_\eta$.
\end{lemma}

\begin{proof}
It is straightforward to verify that $\varrho_\eta$ is indeed a metric on $\cG$ and let us only check that $\varrho_{\rm intr} = \varrho_\eta$. 
First, observe that for any two points $x,y$ on $\cG$ and every path $\cP$ from $x$ to $y$, the following estimate
\begin{align}
|f(x) - f(y)| \le \int_\cP |\nabla f|\rD x \le \int_\cP \sqrt{\frac{\mu}{\nu}}\, \rD x = |\cP|_\eta
\end{align}
holds true and hence $\varrho_{\rm intr} \le \varrho_\eta$. On the other hand, define $f\in H^1_{\loc}(\cG)$ by fixing some $y\in\cG$ and then set $f(x) = \varrho_\eta(x,y)$. It is immediate to see that $f$ is edgewise absolutely continuous and, moreover, $|\nabla f| = \sqrt{\frac{\mu}{\nu}}$ a.e. on $\cG$. Thus $f\in\wh{\cD}_{\loc}$. Moreover, for each $x\in \cG$ we clearly have $\varrho_\eta(x,y) = f(x) - f(y) = f(x)$.
\end{proof}

\begin{remark}
It is easy to see that $\varrho_\eta$ generates the topology of the metric graph $\cG$. Notice also that in the case $\mu=\nu$, $\eta$ coincides with the Lebesgue measure and hence in this case $\varrho_\eta$ is nothing but the length metric $\varrho_0$ on $\cG$. 
\end{remark}

\subsection{Schr\"odinger operators on metric graphs}\label{ss:IV.05}

Of course it makes sense to consider Schr\"odinger-type operators on metric graphs, i.e., operators $-\Delta + q$ (these are widely known as {\em quantum graphs}). Indeed, it is not difficult to replace the divergence form Sturm--Liouville differential expression in \eqref{eq:Hemax} by the general three term Sturm--Liouville expression
\begin{align}\label{eq:SLgenE}
\tau_e = \frac{1}{\mu_e}\Big( - \frac{\rD }{\rD x}\nu_e\frac{\rD}{\rD x} + q_e\Big).
\end{align}
Of course, then one must be careful with the regularity assumptions on the corresponding coefficients. The standard Sturm--Liouville theory deals with the case $\mu_e,1/\nu_e,q_e \in L^1(e)$ together with the positivity of the weights $\mu_e$ and $\nu_e$ (see, e.g., \cite{wei}). In fact, one can go far beyond these assumptions (see, e.g., \cite{sash03}, \cite{egnt}). 

If $q\in L^2_{\loc}(\cG)$ is real valued, then under the assumptions of Hypothesis~\ref{hyp:munu}, $q$ can be considered as an edgewise strongly bounded perturbation w.r.t. $\rH_{e,\max}$ and hence adding such a potential would not change the definition of both the maximal operator $\bH_{\max}$ and the corresponding operators $\bH$ and $\bH_\alpha$. 
However, one can in fact allow $q\in H^{-1}_{\loc}(\cG)$, where $H^{-1}_{\loc}(\cG)$ is the topological dual of $H^1_c(\cG)$ (w.r.t. the usual $L^2$ pairing). 
 In this case, the operator $\rH_{e,\max}$ is defined as follows. First, take a real valued $Q\in L^2_{\loc}(\cG)$ such that $Q' = q$ in the sense of distributions, i.e., 
\begin{align*}
q(f) = -\int_\cG Q(x) \nabla f(x)\rD x 
\end{align*} 
for all $f\in H^1_c(\cG)$. Then on each $e\in\cE$ the differential expression \eqref{eq:SLgenE} can be written in the form
\begin{align}\label{eq:SLdistrE}
\tau_e f = \frac{1}{\mu_e}\Big( - (f^{[1]})'  - \frac{Q_e}{\nu_e} f^{[1]} - \frac{Q_e^2}{\nu_e}f\Big), \qquad f^{[1]}:= \nu_e f' - Q_e f,
\end{align}
where $f^{[1]}$ is often called a {\em quasi-derivative} of $f$. The maximal operator $\rH_{e,\max}$ is defined in $L^2(e;\mu)$ by the above differential expression on the domain 
\[
\dom(\rH_{e,\max}) = \{f\in L^2(e;\mu)|\ f,f^{[1]}\in AC(e),\ \tau_e f\in L^2(e;\mu)\}.
\]
This operator is closed and well-defined (see, e.g., \cite{sash03}, \cite{hrmy}) and hence one can further proceed as in Section~\ref{ss:IV.02}. Indeed, the above definition implies that the quantities \eqref{eq:tr_fe} together with  
 \begin{align}\label{eq:tr_f^1e}
\partial_{\vec{e}}^{[1]} f (\Ei) & := \lim_{x_e \to \Ei} f^{[1]}(x_e), & 
\partial_{\vec{e}}^{[1]} f (\Et) & := \lim_{x_e \to \Et} f^{[1]}(x_e),
 \end{align}
 are well defined, and hence we can introduce Kirchhoff conditions at the vertices
 \begin{align}\label{eq:kirchhoff-quasi}
\begin{cases} f\ \text{is continuous at}\ v,\\[1mm] 
\sum\limits_{\vec{e}\in \vec{\cE}_v} \pi_v({\vec{e}})\, \partial_{\vec{e}}^{[1]} f(v) = 0, \end{cases}\qquad v\in \cV.
\end{align} 
Let us emphasize that the quasi-derivatives and hence $\partial_{\vec{e}}^{[1]} f $ do depend on the choice of the primitive $Q$ of $q$, but the Kirchhoff conditions~\eqref{eq:kirchhoff-quasi} are easily seen to be independent of this choice. The remaining details follow analogously to the considerations of Section~\ref{ss:IV.02} and we leave them  to the interested reader. In the following we shall denote the corresponding differential expression by $-\Delta + q$ and the pre-minimal, minimal and maximal operators by $\bH_q'$, $\bH_q^0$ and, respectively, $\bH_q$.

\begin{remark}[$\delta$-couplings as $\delta$-potentials]\label{rem:LaplDeltaII}
Let us demonstrate how $\delta$-couplings \eqref{eq:vert-alpha} can be considered as a singular potential. For a given model $(\cV,\cE,|\cdot|)$ of $\cG$ consider an edgewise constant function $Q = Q_\alpha$ on $\cG$ such that 
\begin{equation} \label{eq:BoundaryOperator}
\sum_{\vec{e}\in\vec{\cE}_v} \pi_v({\vec{e}})\,Q_\alpha(e) = \alpha(v) 
\end{equation}
for all $v\in\cV$\footnote{Notice that such a function $Q_\alpha$ indeed exists for every $\alpha\colon \cV \to \R$. Namely, after fixing an orientation on the edge set $\cE$, we can identify $\R^\cE$ and $\R^\cV$ with the spaces of $1$-chains and $0$-chains on the infinite, combinatorial graph $\cG_d = (\cV, \cE)$  (see \cite[Chap.~I.2]{soardi}). Moreover, \eqref{eq:BoundaryOperator} defines an operator $\partial \colon \R^\cE \to \R^\cV$, which is nothing but the usual boundary operator from $1$-chains to $0$-chains in this context. The surjectivity of this operator on infinite, locally finite graphs is well-known and can be understood as the corresponding $0$-th homology group being trivial (see, e.g., \cite[Prop.~11.1.3]{geog} or \cite[p.9]{soardi}).}. 
Here we identified $Q_\alpha$ with the function $Q_\alpha\colon \cE\to\R$, and $\pi_v$ is the orientation function on $\vec{\cE}_v$. 
It is straightforward to see that
\begin{align*}
q_\alpha(f) := &- \int_\cG Q_\alpha(x)\nabla f(x)\rD x = -\sum_{e\in\cE} Q_\alpha(e)(f(\Et) - f(\Ei)) \\
&=\sum_{v\in\cV} f(v)\sum_{\vec{e}\in\vec{\cE}_v} \pi_v(\vec{e})Q_\alpha(e) = \sum_{v\in\cV} \alpha(v)f(v)
\end{align*}
for all $f\in H^1_c(\cG)$ and hence $q_\alpha = Q_\alpha' = \sum_{v}\alpha(v)\delta_v$ in the sense of distributions. For this choice of $Q_\alpha$, the differential expression \eqref{eq:SLdistrE} clearly coincides with \eqref{eq:Hemax} (since $Q_e$ is constant on $e$). The last condition in~\eqref{eq:kirchhoff-quasi} reads 
\begin{align*}
0 = \sum_{\vec{e}\in\vec{\cE}_v}  \pi_v(\vec{e}) \partial_{\vec{e}}^{[1]} f(v)
& = \sum_{\vec{e}\in\vec{\cE}_v} \pi_v(\vec{e}) \big(\nu_{\vec{e}}(v) \pi_v(\vec{e}) \partial_{\vec{e}}f (v) - Q(e) f(v)\big)\\
&= \sum\limits_{\vec{e}\in \vec{\cE}_v} \nu_{\vec{e}}(v)\partial_{\vec{e}} f(v) 
- f(v)\sum\limits_{\vec{e}\in \vec{\cE}_v} \pi_v(\vec{e}) Q(e) \\
&= \sum\limits_{\vec{e}\in \vec{\cE}_v} \nu_{\vec{e}}(v)\partial_{\vec{e}} f(v) 
- \alpha(v)f(v),
\end{align*}
which is nothing but the last condition in \eqref{eq:vert-alpha}.  Therefore, the Laplacian with boundary conditions~\eqref{eq:vert-alpha} can be written  as (see also \eqref{eq:QFalpha})
\begin{align}\label{eq:LaplAlpha}
-\Delta + \sum_{v\in\cV}\alpha(v)\delta_v.
\end{align}
\end{remark}

\section{The Glazman--Povzner--Wienholtz theorem on metric graphs}\label{sec:GPWmetr}

Our main result in this paper is the following statement:

\begin{theorem}\label{th:GPWcont}
Let $(\cG,\mu,\nu)$ be a weighted metric graph satisfying Hypothesis~\ref{hyp:munu} and such that $(\cG,\varrho_\eta)$ is complete. Assume that $q\in H^{-1}_{\loc}(\cG)$ is a real distribution such that the minimal operator $\bH^0_q$ associated with $-\Delta + q$ in $L^2(\cG;\mu)$ is bounded from below. Then $\bH^0_q$ is self-adjoint.
\end{theorem}

As an immediate corollary we arrive at the Carleman--Friedrichs result for Kirchhoff Laplacians. For simplicity we state it for $q\in L^1_{\loc}(\cG)$.

\begin{corollary}\label{cor:GaffneyCont}
Let $(\cG,\mu,\nu)$ be a weighted metric graph satisfying Hypothesis~\ref{hyp:munu} and such that $(\cG,\varrho_\eta)$ is complete. If $q\in L^1_{\loc}(\cG)$ is bounded from below, $q\ge -c$ a.e. on $\cG$ for some $c>0$, then the maximal operator $\bH_q = -\Delta +q$ is self-adjoint and coincides with the minimal operator $\bH_q^0$. 
\end{corollary}

\begin{proof}
$\bH^0$ is nonnegative and hence Theorem~\ref{th:GPWcont} can be applied.
\end{proof}

\begin{remark}\label{rem:StuL2}
By employing the version of Yau's~$L^2$-Liouville theorem for strongly local Dirichlet forms from~\cite{stu}, one can prove Corollary~\ref{cor:GaffneyCont} assuming only that the weights $\mu$ and $1/\nu$ are $L^1_{\loc}$. In this respect, we refer to~\cite[Chap.~7.1]{kn21} where the case of edgewise constant weights and $q\equiv 0$ is considered. Other proofs of Corollary~\ref{cor:GaffneyCont} in this special case can be found in \cite[Cor.~4.9]{ekmn} and \cite[Theorem~3.49]{hae}.
\end{remark}

\begin{proof}[Proof of Theorem \ref{th:GPWcont}]
To simplify the exposition, we restrict ourselves to the case of $\delta$-potentials, that is, we provide the proof in the case $q = \sum_v \alpha(v)\delta_v$ (see Remark~\ref{rem:LaplDeltaII}). The case of a general real $q\in H^{-1}_{\loc}(\cG)$ follows line to line the case of $\delta$-potentials, however, some notations become a bit more cumbersome.

Assume that $\bH^0_\alpha$ is bounded from below, that is, there is $c\ge 0$ such that
\begin{align}\label{eq:HaSBD}
\big\langle\bH^0_\alpha f,f \big\rangle_{L^2(\cG;\mu)} \ge -c\|f\|_{L^2(\cG;\mu)}^2
\end{align}
for all $f\in \dom(\bH^0_\alpha)$. If $\bH^0_\alpha$ is not self-adjoint, then $\ker(\bH_\alpha + (c+1)\rI)\neq \{0\}$, that is, there exists $0\neq u\in \dom(\bH_\alpha)$ such that $\bH_\alpha u + (c+1)u = 0$.  

Next, fix a model of $(\cG,\mu,\nu)$. Clearly,  we can assume that our model is simple (no loops or multiple edges). 
Let $(\cG_n)_{n\ge 0}$ be a compact exhaustion of $\cG$, that is, $\cG_n$ is a finite connected subgraph of $\cG$, $\cG_n\subsetneq\cG_{n+1}$ for each $n$ and $\cup_{n\ge 0}\cG_n = \cG$. Suppose that $\varphi_n \in L^2_c(\cG;\mu)$ is real-valued, continuous and such that $\varphi_n|_{\cG_n} \equiv 1$  and $u\varphi_n\in\dom(\bH^0_\alpha)$. 
Then using \eqref{eq:HaSBD} we get
\begin{align}\label{eq:contradiction}
\|u\|^2_{L^2(\cG_n;\mu)} \le \|u\varphi_n\|^2_{L^2(\cG;\mu)} \le \big\langle (\bH_{\alpha}^0 + (c+1)\rI)(u\varphi_n),u\varphi_n\big\rangle_{L^2(\cG;\mu)}. 
\end{align}
On the other hand, on each edge $e\in \cE$ we have
by employing a simple integration by parts
\begin{align*}
\int_e\big(-\Delta (u\varphi_n) + &(c+1)u\varphi_n\big)u\varphi_n\,\mu(\rD x) \\ 
&  = \int_e\Big( -\varphi_n \Delta u - u\Delta \varphi_n - \frac{2}{\mu}\big(\nu u'\varphi_n' \big) + (c+1) u\varphi_n\Big)u\varphi_n\,\mu(\rD x) \\
& = - 2\int_e uu'\, \varphi_n\varphi_n'\,\nu(\rD x) - \int_e  u^2 \varphi_n (\nu\varphi_n')'\,\rD x \\
& = - \frac{1}{2}\int_e  (u^2)' (\varphi_n^2)'\,\nu(\rD x) - \int_e  u^2 \varphi_n (\nu\varphi_n')'\,\rD x \\
& = -  \nu\varphi_n'\, \varphi_nu^2\Big|_{\Ei}^{\Et} + \int_e  \frac{1}{2}u^2 \big(\nu(\varphi_n^2)'\big)' - u^2 \varphi_n (\nu\varphi_n')'\,\rD x \\
& = -  \nu\varphi_n'\, \varphi_nu^2\Big|_{\Ei}^{\Et} + \int_e  u^2  (\varphi_n')^2\,\nu(\rD x) \\
& = - \nu\varphi_n'\, \varphi_nu^2\Big|_{\Ei}^{\Et} + \int_e   \Big(\frac{u\varphi_n'}{\eta}\Big)^2 \mu(\rD x).
\end{align*} 
Summing up over all $e\in\cE$, this together with \eqref{eq:contradiction} lead to the following inequality upon choosing a suitable $\varphi_n$ to neglect the boundary terms
\begin{align}\label{eq:VII.1.02}
\|u\|^2_{L^2(\cG_n;\mu)} \le \Big\|u\,\frac{\nabla \varphi_n}{\eta}\Big\|^2_{L^2(\cG\setminus \cG_n;\mu)}.
\end{align}
Thus, if in addition a sequence of compactly supported functions $(\varphi_n)$ satisfies $|\nabla \varphi_n|\le C\,\eta$ a.e. on $\cG$ for all $n\ge 0$ and with some constant $C>0$ independent of $n$, we arrive at contradiction by sending $n$ to infinity in \eqref{eq:VII.1.02}.

It remains to show that a sequence $(\varphi_n)$ with all the above properties indeed exists. First of all, notice that $u\varphi_n$ is continuous if so are both $u$ and $\varphi_n$. Thus, $u\varphi_n\in\dom(\bH^0_\alpha)$ if $\varphi_n$ is $H^2$ edgewise and satisfies (recall the simplicity of the model assumption)
\begin{align}\label{eq:VII.1.03}
\sum\limits_{e\in \cE_v}\nu_e(v) \partial_{e} \varphi_n(v) = 0
\end{align}
at each vertex $v\in\cV$. 
Since $\cG_n$ is a finite subgraph of a locally finite graph $\cG$, its boundary $\partial \cG_n$ is a finite set. Let us now define $\wt\cG_n\supset \cG_n$ as the closed radius $1$ distance ball of $\cG_n$, that is, $\wt\cG_n = B_1(\cG_n;\eta)$, where 
\[
B_r(\cG_n; \eta) := \{x\in \cG\,|\ \varrho_\eta(x,\cG_n)\le r\},\qquad r>0.
\] 
Taking into account that $\cG$ is complete w.r.t. $\varrho_\eta$, $\wt\cG_n$ is a compact subset of $\cG$ by the Hopf--Rinow theorem for length spaces (see \cite[Theorem~2.5.28]{bbi}). This in particular implies that $\wt\cG_n \cap \cV$ is a finite set. 

Next, consider the function
\begin{align}\label{eq:testfunc}
\wt\varphi_n(x) = \max (0, 1-\varrho_\eta(x,\cG_n)),\qquad x\in \cG.
\end{align}
Clearly, $\wt\varphi_n$ is continuous on $\cG$. Moreover, 
\begin{align}\label{eq:testfuncDiff}
|\nabla \wt\varphi_n| = \eta\ \ \text{a.e. on}\ \ \wt\cG_n\setminus \cG_n
\end{align}
and $\nabla\wt\varphi_n = 0$ inside $\cG_n$ and $\cG\setminus \wt\cG_n$. However, it may happen that $\nabla \wt\varphi_n$ does not satisfy \eqref{eq:VII.1.03} on the vertices inside $\wt\cG_n$. Moreover, $\nabla \wt\varphi_n $ may be discontinuous at the boundaries of $\cG_n$ and $\wt\cG_n$ as well as at the points inside edges in $\wt\cG_n\setminus \cG_n$, where it changes sign. However, it is easy to see that for each edge $e\in \cG$ this may happen at most finitely many times. Applying the Hopf--Rinow theorem once again, one easily concludes that the set of all these points is finite. By refining our chosen model, we can assume without loss of generality that all these points are vertices of $\cG$. Moreover, taking into account that $\cG$ is locally finite, we can use a ``smoothing" (similar to changing the modulus function into a smooth function at $0$) in order to turn $\wt\varphi_n$ into $\varphi_n$, satisfying all the smoothness requirements and keeping the assumption $|\nabla\varphi_n| \le C\eta$ on $\cG$ with some $C>0$ uniform in $n$. 
More precisely, fix a smooth function $\psi \colon [0,1] \to [0,\infty)$ supported on $[0, 1/2]$ with $\psi'(0) =1$ and $\psi(0)=0$. Defining 
\[
	\psi_{v,e}  :=  |e| \psi \Big (\frac{|\cdot-v|}{|e|} \Big),
\] 
for every $v \in \cV$ and $e \in \cE_v$, we can modify $\wt\varphi_n$ in a vicinity of all the vertices of $\wt\cG_n$ by setting on the corresponding edges $e = e_{u,v}$
 \[
	\varphi_n \big|_u^v := \wt\varphi_n \big|_u^v - (\partial_e \wt\varphi_n(v))  \psi_{v,e} - (\partial_e \wt\varphi_n(u))  \psi_{u,e}.
\]
Then by construction $|\nabla \varphi_n| \le (1+C) \eta$ a.e. on $\cG$ with $C := \| \psi'\|_\infty$ and, moreover, $\partial_e\varphi(v) = 0$ for all $e\in\cE_v$ and $v\in \cV$.
 This finishes the proof.
\end{proof}

\begin{remark}\label{rem:GPWmetric}
A few remarks are in order.
\begin{itemize}
\item[(i)]
In the special case of a path graph, that is for Schr\"odinger operators in $L^2(\R)$, Theorem~\ref{th:GPWcont} was proved by different approaches in \cite[Theorem 1 and Rem.~III.2]{akm10} and \cite{hrmy}. Our proof follows the ideas of \cite{akm10}. It might be interesting to check whether the proof of \cite{hrmy}, which is based on the use of supersymmetry, extends to the graph setting.  
\item[(ii)]
It seems the claim of Theorem~\ref{th:GPWcont} can be extended to edge weights $\mu$ and $\nu$ assuming only that $\mu,1/\nu \in L^1_{\loc}(\cG)$. However, the main difficulty is to construct suitable test functions $\varphi_n$ since our construction breaks down at several places. This issue will be addressed elsewhere.
\item[(iii)]
It is tempting to replace the completeness w.r.t $\varrho_\eta$ in Theorem \ref{th:GPWcont} by the one w.r.t. the star path metric $\varrho_m$. In particular, we are convinced that the following analog of Corollary~\ref{cor:GaffneyCont} takes place (cf. \cite[Theorem~4.11]{ekmn} and \cite[Theorem~7.7]{kn21}): {\em if $q\in L^1_{\loc}(\cG)$ is bounded from below and for a finite size model $(\cV,\cE,|\cdot|,\mu,\nu)$ the metric space $(\cV,\varrho_m)$ is complete, where $\varrho_m$ is the path metric \eqref{eq:rho_m} and $m = m_\mu$ is the vertex weight defined by}
\[
m_\mu(v) = \mu(\cE_v) = \int_{\cE_v}\mu(\rD x),\qquad v\in\cV,
\]
{\em then $\bH_q = -\Delta+q$ is self-adjoint}. On the other hand, simple examples show that the analog of Theorem~\ref{th:GPWcont} is not true under this completeness assumption (see Remark~\ref{rem:wouk}(iv)).
\end{itemize}
\end{remark}

\section{The Glazman--Povzner--Wienholtz theorem on weighted graphs}\label{sec:GPWgraph}
As an immediate application of Theorem \ref{th:GPWcont} and the results connecting metric graphs with weighted graphs (see Appendix~\ref{app:WGvsWMG}), we arrive at the following version of the Glazman--Povzner--Wienholtz theorem on graphs.

\begin{theorem}[The Glazman--Povzner--Wienholtz theorem on graphs]\label{th:GPWdiscr}
Let $b$ be a locally finite graph over $(\cV,m)$ and let $\varrho$ be an intrinsic metric which generates the discrete topology on $\cV$.  
Assume also that $\alpha\colon\cV\to \R$ is such that the minimal Schr\"odinger operator $\rh_\alpha^0$ associated with \eqref{eq:SchrDiscr} is bounded from below in $\ell^2(\cV;m)$. 
If $(\cV,\varrho)$ is complete as a metric space, then $\rh^0_\alpha$ is self-adjoint and $\rh^0_\alpha=\rh_\alpha$.
\end{theorem}

\begin{proof}
We prove the claim in three steps. The strategy is to show first that it suffices to consider intrinsic path metrics having finite jump size and then to exploit the connection between Laplacians on graphs and metric graphs established in \cite{ekmn} and also in \cite{kn21}.

(i) Let $\varrho$ be an intrinsic metric which generates the discrete topology on $\cV$ and such that $(\cV, \varrho)$ is complete. 
Consider the weight $p \colon \cV  \times \cV \to [0, \infty)$ given by 
\begin{align*}
p(x,y) := \begin{cases} \varrho(x,y), & b(x,y)>0\\ 0, & b(x,y)=0\end{cases}.
\end{align*}
By construction, $p$ is an intrinsic weight w.r.t. $(\cV,m;b)$ and the associated strongly intrinsic path metric $\wt{\varrho} = \varrho_p$ satisfies $ \rho \le \wt{\varrho}$.
Moreover, using the fact that both $\wt{\varrho}$ and $\varrho$ generate the discrete topology on $\cV$, the completeness of $(\cV, \wt{\varrho})$ follows by comparison. 

(ii) Suppose now that  $\varrho = \varrho_p$ is a general intrinsic path metric with weight function $p \ge 0$ such that $(\cV, \varrho)$ is complete.
 By the discrete Hopf--Rinow theorem (see \cite[Theorem~A.1]{hkmw13}, \cite[Theorem~3.4]{kel15}), 
 the completeness is equivalent to the fact that 
\begin{equation} \label{eq:gaffney_proof2} 
\sum_{n\ge 0} p(v_n, v_{n+1}) =  \infty
\end{equation}
for any infinite path $\cP = (v_0, v_1, v_2,\dots)$ (i.e., $b(v_n, v_{n+1}) >0$ for all $n\ge 0$, see \eqref{eq:pathmetric}). However, introducing the new weight function $\wt p := \min\{ 1, p \}$, we arrive at another path metric $\wt \varrho := \varrho_{\tilde{p}}$, which is strongly intrinsic with respect to $(\cV,m;b)$ (by construction) and, moreover, has jump size at most $1$, that is, 
\begin{align}\label{eq:jumpsizes}
s(\varrho):= \sup\{\varrho(u,v)\,|\, u,v\in\cV, b(u,v)>0\}  \le 1.
\end{align}
 It is not hard to show (e.g., by employing the Hopf--Rinow theorem once again)  that $(\cV, \varrho)$ is complete exactly when so is $(\cV, \wt \varrho)$.
 
(iii) Thus, it remains to prove the claim assuming that $\varrho$ is an intrinsic path metric of finite jump size, $s(\varrho)<\infty$. First of all, by \cite[Lemma~6.27]{kn21} (see also Theorem~\ref{th:metricDvsC}), for a locally finite graph $b$ over $(\cV,m)$ there exists a metric graph $\cG = (\cV,\cE,|\cdot|)$ equipped with two weights $\mu,\nu\colon\cE\to (0,\infty)$ such that 
\[
\sup_{e\in\cE}\int_e \eta(\rD x) = \sup_{e\in\cE}\eta(e) = \sup_{e\in\cE} |e|\sqrt{\frac{\mu(e)}{\nu(e)}} <\infty,
\]
and the weights $m$ and $b$ can be represented by \eqref{eq:Mbndry} and \eqref{eq:b_weight}. Moreover, the metric $\varrho$ can be obtained as the path metric induced by the intrinsic path metric $\varrho_\eta$, that is, $\varrho = \varrho_\eta|_{\cV\times\cV}$. 

Consider now the Laplacian $\bH^0_\alpha$ on the weighted metric graph $(\cG,\mu,\nu)$. It turns out that the properties of $\bH^0_\alpha$ and $\rh^0_\alpha$ are closely connected. First of all, by Theorem~\ref{th:BTmain}(ii), 
 the operator $\bH^0_\alpha$ is bounded from below since so is $\rh^0_\alpha$. Moreover, $(\cG,\varrho_\eta)$ is complete since so is $(\cV,\varrho)$ (this follows by using the Hopf--Rinow theorem and comparing the geodesic completeness on $(\cG,\varrho_\eta)$ and $(\cV,\varrho)$). Now, applying the Glazman--Povzner--Wienholtz theorem for metric graphs, we conclude that  $\bH^0_\alpha$ is self-adjoint. However, by Theorem~\ref{th:BTmain}(i), 
 the operators $\bH_\alpha^0$ and $\rh^0_\alpha$ are self-adjoint only simultaneously.
\end{proof}

\begin{remark}
To the best of our knowledge the Glazman--Povzner--Wienholtz theorem for graphs was established first in \cite[Theorem~1.3]{mil11} and \cite[Theorem~6.1]{toha10} (however, under the additional bounded geometry assumption on $(\cV,b)$, that is, if $\sup_{v\in\cV}\deg(v)<\infty$, and for some particular choices of path metrics), and then in \cite[Theorem~2.16]{gks15}, which allows non-locally finite graphs (see also \cite{schm20}). 
\end{remark}

Let us finish this section by applying Theorem~\ref{th:GPWdiscr} in the case when the vertex weight $m$ is constant on $\cV$, i.e., $m\equiv \id$. Let us rewrite \eqref{eq:SchrDiscr} in this case as
\be\label{eq:jacopgraph}
 (\tau f) (v) = a (v) f(v) - \sum_{u \in \cV} b(u,v) f(u),\qquad v\in\cV. 
\ee 
Here, as before, $b$ is a locally finite graph over $\cV$. Notice that \eqref{eq:jacopgraph} is nothing but \eqref{eq:SchrDiscr} with $m\equiv \id$ and $a(v) = \alpha(v) + \sum_{u \in \cV} b(u,v)$, $v\in\cV$. Let us denote the corresponding minimal and maximal operators associated with \eqref{eq:jacopgraph} in $\ell^2(\cV)$ by $\rJ^0 = \rJ^0_{a,b}$ and $\rJ = \rJ_{a,b}$.

\begin{remark}[Jacobi matrices on graphs]\label{rem:jacop}
Sometimes the second order difference expression \eqref{eq:jacopgraph} is called a {\em Jacobi matrix on a graph} (cf., e.g., \cite{ady20}, \cite{abs20}). Indeed, if $\cV = \Z_{\ge 0}$ and $b$ is a path graph over $\Z_{\ge 0}$, that is, $b(n,m) >0$ exactly when $|n-m| = 1$, then \eqref{eq:jacopgraph} can be written as a Jacobi (tri-diagonal) matrix
\begin{align}\label{eq:1Djacop}
 \begin{pmatrix} a_0 & -b_0 & 0 & 0 & \dots \\
-b_0 & a_1 & -b_1 &  0 & \dots \\
0& -b_1 & a_2 & -b_2 & \dots \\
 0 & 0 & -b_2 & a_3 & \dots \\
\dots & \dots & \dots & \dots & \dots \\
\end{pmatrix},
\end{align}
where we set $a_n = a(n)$ and $b_n = b(n,n+1)$, $n\in\Z_{\ge 0}$.
\end{remark}

With \eqref{eq:jacopgraph}, let us associate the path metric $\varrho_b$ defined by the weights
\begin{align}\label{eq:WeightJacop1}
p_b(u,v) = \begin{cases} \frac{1}{\sqrt{b(u,v)\,\max(\deg(u),\deg(v))}}, & b(u,v)\neq 0,\\[2mm] 0, & b(u,v) = 0.\end{cases}
\end{align}
Here $\deg(v) := \#\{ u \in \cV\,| \, b(u,v) >0 \}$ for every $v \in \cV$, which coincides with the degree of $v$ in the simple graph $\cG_b$ associated to $b$ (see Remark~\ref{rem:simplevsmult}).

\begin{theorem}\label{th:WoukGraph}
If the minimal operator $\rJ^0$ is bounded from below and $(\cV,\varrho_b)$ is complete as a metric space,
then the operator $\rJ^0$ is self-adjoint in $\ell^2(\cV)$ and $\rJ^0 = \rJ$.
\end{theorem}

\begin{proof}
The proof is  a straightforward application of Theorem~\ref{th:GPWdiscr}. Indeed, we only need to show that the metric $\varrho_b$ is intrinsic. Recalling that $m\equiv \id$, we get
\begin{align*}
\sum_{u\sim v} b(u,v)p_b(u,v)^2 = \sum_{u\sim v} \frac{1}{\max(\deg(u),\deg(v))} \le \sum_{u\sim v} \frac{1}{\deg(v)} = 1, 
\end{align*}
for all $v\in\cV$. 
\end{proof}

\begin{remark}\label{rem:wouk}
A few remarks are in order.
\begin{itemize}
\item[(i)]
If the graph $b$ has bounded geometry, that is, $\sup_\cV \deg(v)<\infty$, then \eqref{eq:WeightJacop1} can be replaced by the equivalent path metric defined by the weights
\begin{align}\label{eq:WeightJacop2}
\wt{p}_b(u,v) = \frac{1}{\sqrt{b(u,v)}}, \quad u\sim v.
\end{align}
\item[(ii)]
Under the bounded geometry assumption, Theorem~\ref{th:WoukGraph} was proved earlier in \cite[Theorem~1.3]{mil11} and \cite[Theorem~6.1]{toha10}. 
\item[(iii)]
Let us also mention that Theorem~\ref{th:WoukGraph} can be seen as an extension of the self-adjointness test for Jacobi matrices \eqref{eq:1Djacop} established by A.~Wouk in~\cite{wouk} (see also~\cite[Problem~I.4]{Akh}). Indeed, for \eqref{eq:1Djacop}, the completeness of $\Z_{\ge 0}$ with respect to $\varrho_b$ is equivalent to the following condition
\begin{align}\label{eq:woukJacop}
\sum_{n\ge 0}\frac{1}{\sqrt{b_k}} = \infty.
\end{align}
\item[(iv)]
 It is well-known that one cannot remove \eqref{eq:woukJacop} in the Wouk test, that is, semiboundedness of \eqref{eq:1Djacop} does not guarantee its essential self-adjointness\footnote{For instance, choose the Jacobi parameters in \eqref{eq:1Djacop} in the form $a_n = m_n(l_{n-1} + l_n)$, $b_n = l_n\sqrt{m_n m_{n+1}}$ with $m_n = (n+1)^{c_1}$ and $l_n = (n+1)^{c_2}$ for some positive constants $c_1,c_2>0$. It is straightforward to verify that the corresponding minimal operator is nonnegative. However, it is not self-adjoint  in $\ell^2$ if and only if $\sum_{n\ge 0} \frac{1}{m_n}(\sum_{k=0}^{n-1} \frac{1}{l_k})^2 <\infty$ (cf.~\cite[Theorem~0.5]{Akh}). The latter holds exactly when $c_1>1$ and $c_1+2c_2>3$. Notice that \eqref{eq:woukJacop} holds true exactly when $c_1+c_2\le2$.}. Since $m\equiv 1$, the corresponding path metric $\varrho_m$ given by \eqref{eq:rho_m} is  the combinatorial distance and clearly $(\Z_{\ge 0},\varrho_m)$ is complete. Therefore, this shows that one cannot replace the completeness w.r.t. an intrinsic metric  in Theorem~\ref{th:GPWdiscr} by the completeness w.r.t. $\varrho_m$. 
\item[(v)] Using the isometric isomorphism
\[
\begin{array}{cccc}
U_m\colon & \ell^2(\cV;m) & \to & \ell^2(\cV)\\
 & f & \mapsto & \sqrt{m}f 
\end{array}, 
\]
it is straightforward to see that every Schr\"odinger operator $\rh_\alpha$ associated with \eqref{eq:SchrDiscr} is unitarily equivalent to some $\rJ$, that is, $U_m\rh_\alpha U_m^{-1}$ is the maximal operator generated by a Jacobi matrix  in $\ell^2(\cV)$.
\end{itemize}
\end{remark}

\begin{example}\label{ex:CombLaplPert}
We would like to conclude by considering perturbations of the combinatorial Laplacian~\eqref{eq:LaplComb} considered in Example~\ref{ex:CombLapl}. Namely, let 
\begin{align}
(L_\alpha f)(v) = \alpha(v)f(v) + \sum_{u\sim v} f(v) - f(u) = (\deg(v) + \alpha(v)) f(v) - \sum_{u\sim v}f(u),
\end{align}
for some $\alpha\colon \cV\to \R$. By the discrete version of the Carleman--Friedrichs theorem (see \cite[Theorem~5]{kl12}), the corresponding maximal Schr\"odinger operator $\rh_\alpha$ is self-adjoint whenever $\alpha$ is bounded from below, $\inf_{v\in\cV}\alpha(v) > -\infty$. Theorem~\ref{th:WoukGraph} enables us to replace the latter condition by the semiboundedness of the minimal operator $\rh_\alpha^0$ if $(\cV,\varrho_b)$ is complete. Notice that in this example $b$ is nothing but the adjacency matrix and hence the weight \eqref{eq:WeightJacop1} is simply given by  
\begin{align}\label{eq:WeightJacopComb}
p_b(u,v) = \frac{1}{\sqrt{\max(\deg(u),\deg(v))}}
\end{align}
when $u\sim v$, and $p_b(u,v)=0$ otherwise. It turns out that it is impossible to remove this completeness assumption. Indeed, this fact readily follows from Remark~\ref{rem:wouk}(iv). Namely, consider an {\em antitree}\footnote{Consider a simple graph $\cG_d = (\cV, \cE)$. Fix a root vertex $o \in \cV$ and then order the graph with respect to the combinatorial spheres\label{not:combsphere} $S_n$, $n \in \Z_{\ge 0}$ ($S_n$ consists of all vertices $v\in\cV$ such that the combinatorial distance from $v$ to the root $o$, that is, the combinatorial length of the shortest path connecting $v$ with $o$, equals $n$; notice that $S_0=\{o\}$). A connected simple rooted  graph $\cG_d$ is called an {\em antitree}\label{not:antitree2} if every vertex in the combinatorial sphere $S_n$, $n\ge 1$, is connected to all  vertices in $S_{n-1}$ and $S_{n+1}$ and no vertices in $S_k$ for all $|k-n|\neq 1$. Clearly, every antitree is uniquely determined by its sphere numbers $(s_n)_{n\ge 0}$, $s_n = \#S_n$.} 
with sphere numbers $(s_n)_{n\ge 0} \subseteq \Z_{\ge 1}$ and assume that $\alpha\colon \cV\to\R$ is radially symmetric, that is, $\alpha|_{S_n} \equiv \alpha_n$ for all $n\ge 0$. Using \cite[Theorem~4.1]{brke13}, one can easily show that the Schr\"odinger operator $\rh_\alpha^0$ is self-adjoint/bounded from below exactly when so is the Jacobi matrix \eqref{eq:1Djacop} with the Jacobi parameters
\begin{align}\label{eq:JacopAT}
b_n & = \sqrt{s_ns_{n+1}}, & a_n & = \alpha_n + s_{n-1} + s_{n+1}.
\end{align} 
Now let $m_n = \floor{(n+1)^{c_1}}$ and $l_n = \sqrt{\gamma_{n+1}\gamma_n}$, where $\gamma_n = \floor{(n+1)^{c_2}}$ for all $n\in\Z_{\ge 0}$ with some positive constants $c_1,c_2>0$. Here $\floor{\cdot}$ is the floor function. Next, setting 
\begin{align*}
s_{n} & = m_n \gamma_n = \floor{(n+1)^{c_1}}\floor{(n+1)^{c_2}}, & \alpha_n & = m_{n}(l_{n-1} + l_n) - (s_{n-1} + s_{n+1}),
\end{align*} 
for all $n\in\Z_{\ge 0}$ with $s_{-1}=0$, the corresponding Jacobi parameters~\eqref{eq:JacopAT} admit the factorisation $b_n = l_n\sqrt{m_n m_{n+1}}$ and $a_n =  m_{n}(l_{n-1}+l_n)$ as in the footnote example on the previous page. Therefore, we end up with the minimal Schr\"odinger operator, which is always nonnegative. However, if $c_1>1$ and $c_1+2c_2>3$, then it is non-self-adjoint. Taking into account that  $\deg|_{S_n} = s_{n-1}+s_{n+1}$ and $s_n\approx n^{c_1+c_2}$ as $n\to \infty$, this example also demonstrates that the completeness condition (w.r.t. the weight \eqref{eq:WeightJacopComb}) in Theorem~\ref{th:WoukGraph} is sharp in the case of the combinatorial Laplacian.

\end{example}

\appendix

\section{Weighted graphs vs. weighted metric graphs}\label{app:WGvsWMG}

Consider a weighted metric graph $(\cG,\mu,\nu)$ with edgewise constant weights, that is, for some model $(\cV,\cE,|\cdot|)$ of $\cG$ the weights $\mu$ and $\nu$ are constant in the interior of each edge $e\in\cE$. Hence we can identify $\mu$ and $\nu$ with the functions
\[
\mu,\ \nu\colon\cE\to (0,\infty).
\]
We shall also assume that the model has finite intrinsic size, that is, 
\begin{align}\label{eq:finsize}
\sup_{e\in\cE} |e|\sqrt{\frac{\mu(e)}{\nu(e)}} < \infty.
\end{align}
Consider the corresponding Laplacian $\bH_\alpha$ equipped with vertex conditions \eqref{eq:vert-alpha} for some $\alpha\colon\cV\to\R$. Suppose $f$ is an edgewise affine continuous function on $(\cV,\cE,|\cdot|)$ satisfying \eqref{eq:vert-alpha} on $\cV$. Then it satisfies the last equality in \eqref{eq:vert-alpha} if its slopes satisfy 
\begin{align*}
\sum_{u\sim v} \sum_{\vec{e} \in \vec{\cE}_{u}\colon e\in\cE_v } \frac{\nu(e)}{|e|} \big (f(u)- f(v) \big)    = \alpha(v)f(v),
\end{align*}
 (take into account that we allow multigraphs and any function $f$ must be constant on loop edges). Introducing the Hilbert space $\ell^2(\cV;m)$ with the vertex weight
\begin{align}\label{eq:Mbndry}
m\colon v\mapsto \sum_{\vec{e}\in\vec{\cE}_v} |e|\mu(e),\qquad v\in\cV,
\end{align}
the finite size condition \eqref{eq:finsize} ensures that $f\in L^2(\cG;\mu)$ exactly when  $\f := f|_\cV$ belongs to $\ell^2(\cV;m)$ (see~\cite[Section 4.3]{kn21} and \cite[Rem.~3.8]{ekmn}). This indicates a close relationship between the operator $\bH_\alpha$ and the weighted discrete Schr\"odinger operator $\rh_\alpha$ defined in $\ell^2(\cV;m)$ by the difference expression \eqref{eq:SchrDiscr} with $m$ given by \eqref{eq:Mbndry} and the graph $b$ over $(\cV,m)$ defined by 
\begin{align}\label{eq:b_weight}
b(u,v)  
=  \begin{cases}  
                 \sum_{\vec{e} \in \vec{\cE}_{u}\colon e\in\cE_v } \frac{\nu(e)}{|e|}, & u \neq v,  \\[1mm] 
                  \quad 0,  & u = v.
       \end{cases}
\end{align}
These connections are by no means new and here we are going to recall only the facts needed in the proof of Theorem~\ref{th:GPWdiscr}. First, we need the result relating their spectral properties.

\begin{theorem}\label{th:BTmain}
Let $\cG = (\cV,\cE,|\cdot|)$ be a metric graph equipped with two edge weights $\mu,\nu\colon \cE\to(0,\infty)$ such that \eqref{eq:finsize} is satisfied. Let also $\alpha\colon\cV\to\R$ and $\bH^0_\alpha$ be the corresponding minimal operator. If $\rh_\alpha^0$ is the minimal Schr\"odinger operator associated with \eqref{eq:SchrDiscr}, \eqref{eq:Mbndry}--\eqref{eq:b_weight} in $\ell^2(\cV;m)$, then:
\begin{itemize}
\item[(i)] The deficiency indices of $\bH_\alpha^0$ and $\rh_\alpha^0$ are equal and 
\begin{align}\label{eq:def+=-}
\Nr_+(\bH_\alpha^0) = \Nr_-(\bH_\alpha^0) = \Nr_\pm(\rh_\alpha^0)\le \infty.
\end{align}
In particular, $\bH_\alpha^0$ is self-adjoint if and only if $\rh_\alpha^0$ is self-adjoint.
\item[(ii)] 
The operator $\bH_\alpha^0$ is lower semibounded if and only if the operator $\rh^0_\alpha$ is lower semibounded. 
\end{itemize}
\end{theorem} 

\begin{remark}
This result was first proved in \cite{KM10} for 1d Schr\"odinger operators with point interactions and later it was extended to Laplacians on metric graphs \cite[Theorem~3.5]{ekmn} (see also \cite[Chap.~3]{kn21}). The list of equivalences is much longer and we refer to  \cite{ekmn}, \cite{kn21} for further details. Let us also mention that Theorem~\ref{th:BTmain} can be viewed as a part of a general relationship between a self-adjoint (or, more generally, symmetric) extension of an abstract symmetric operator and the boundary operator parameterizing this extension in some appropriate sense (see \cite{DM91} for details). In particular, the operator $\rh_\alpha^0$ indeed serves as a boundary operator of $\bH_\alpha^0$ with respect to a particular boundary triplet constructed in \cite{ekmn},  \cite{KM10}, \cite{kn21}. The general theory of boundary triplets and corresponding Weyl functions is the modern language of extension theory of symmetric operators in Hilbert spaces, which can partially be seen as a far-reaching development of the Birman--Krein--Vishik theory.
\end{remark}

Lemma~\ref{lem:IntrMetWMG} states that the metric $\varrho_\eta$ defined by \eqref{eq:compute_metric} coincides with the (maximal) intrinsic metric $\varrho_{\rm intr}$ of the Dirichlet Laplacian $\bH_D$ on $(\cG,\mu,\nu)$. It is straightforward to verify that the induced path metric on $\cV$,
\[
\varrho_\cV(u,v) := \varrho_\eta(u,v),\qquad u,v\in\cV
\]
is intrinsic w.r.t. $(\cV,m;b)$ (see Section~\ref{ss:III.02}). It turns out that the class of weighted metric graphs with edgewise constant weights is sufficiently large not only to obtain a wide class of Schr\"odinger operators on graphs (in the sense of \eqref{eq:Mbndry}--\eqref{eq:b_weight}), but also intrinsic path metrics. 

\begin{theorem}\label{th:metricDvsC}
Let $b$ be a locally finite, connected graph over $(\cV, m)$ equipped with a strongly intrinsic path metric $\varrho$. Assume also that 
$\varrho$ has finite jump size, 
\begin{align*}
s(\varrho) = \sup\{ \varrho(u,v)\,|\, u,v\in\cV, b(u,v)>0\}  < \infty.
\end{align*} 
Then there exists a metric graph $\cG = (\cV,\cE,|\cdot|)$ equipped with two weights $\mu,\nu\colon\cE\to (0,\infty)$ such that \eqref{eq:finsize} is satisfied, $m$ and $b$ have the form \eqref{eq:Mbndry} and \eqref{eq:b_weight}, respectively, and, moreover, $\varrho$ coincides with the induced path metric $\varrho_\cV = \varrho_\eta|_{\cV\times\cV}$.
\end{theorem}

\begin{remark}
Theorem~\ref{th:metricDvsC} can be found in \cite[Chap.~6]{kn21} (see Lemma~6.27 and Theorem~6.30 there). Notice that upon restricting to weighted metric graphs with equal weights, i.e., $\mu=\nu$, and also assuming some normalization (the so-called {\em canonical models}: no multiple edges and all loops have the same length $1$, see  \cite{kn21}) the correspondence in Theorem~\ref{th:metricDvsC} becomes in a certain sense bijective. Let us also emphasize that some versions of this result have been used earlier in \cite{fo13}, \cite{hua}.
\end{remark} 

\end{document}